	\documentclass[journal,twocolumn, 8pt,english]{IEEEtran}
	
	\usepackage{float}
	\restylefloat{table}
	\usepackage{pseudocode}
	\usepackage{algorithmicx}
	\usepackage{epsfig}
	\usepackage{pdfpages}
	\usepackage{graphicx}
	\usepackage{multirow}
	\usepackage{amsmath,graphicx}
	\usepackage{amssymb}
	\usepackage{epstopdf}
	\usepackage{subfigure}
	\usepackage{color}
	\usepackage{amsthm}
	\newtheorem{thm}{Theorem}

	\theoremstyle{definition}
	
	\theoremstyle{remark}

	\def\be{\begin{equation}}
	\def\bmat{\begin{matrix}}
	\def\emat{\end{matrix}}
	\def\bea{\begin{eqnarray}}
	\def\beas{\begin{eqnarray*}}
	\def\eea{\end{eqnarray}}
	\def\eeas{\end{eqnarray*}}

	\def\bi{\begin{itemize}}
	\def\ee{\end{equation}}
	\def\ei{\end{itemize}}

	\def\z1{z^{-1}}

	\usepackage[hypcap]{caption}

\newcommand{\mbf}[1]{{\mathbf #1}}

\newcommand{\bs}[1]{\boldsymbol #1}

\def\bea{\begin{eqnarray}}
\def\beas{\begin{eqnarray*}}
\def\eea{\end{eqnarray}}
\def\eeas{\end{eqnarray*}}

	\def\be{\begin{equation}}
	\def\bmat{\begin{matrix}}
	\def\emat{\end{matrix}}
	\def\bea{\begin{eqnarray}}
	\def\beas{\begin{eqnarray*}}
	\def\eea{\end{eqnarray}}
	\def\eeas{\end{eqnarray*}}

	\def\bi{\begin{itemize}}
	\def\ee{\end{equation}}
	\def\ei{\end{itemize}}

	\def\z1{z^{-1}}

	\def\bmat{\begin{matrix}}
	\def\emat{\end{matrix}}

\begin{document}

\title{ Subspace based low rank \& joint sparse matrix recovery}
\author{Sampurna Biswas, Sunrita Poddar, Soura Dasgupta, Raghuraman Mudumbai, and Mathews Jacob
\thanks{The authors are with the Department of Electrical and Computer Engineering, Univ. Iowa, IA, USA.  (emails: \{sampurna-biswas,soura-dasgupta,mathews-jacob\}@uiowa.edu).
\newline This work is in part  supported by US NSF grants  EPS-1101284, ECCS-1150801, CNS-1329657, CCF-1302456,  CCF-1116067, NIH 1R21HL109710-01A1, ACS RSG-11-267-01-CCE,  and ONR grant N00014-13-1-0202.}}
\maketitle

	\textit{Abstract:} \textbf{We consider the recovery of a low rank and jointly sparse matrix from under sampled measurements of its columns. This problem is highly relevant in the recovery of dynamic MRI data with high spatio-temporal resolution, where each column of the matrix corresponds to a frame in the image time series; the matrix is highly low-rank since the frames are highly correlated. Similarly the non-zero locations of the matrix in appropriate transform/frame domains (e.g. wavelet, gradient) are roughly the same in different frame. The superset of the support can be safely assumed to be jointly sparse. Unlike the classical multiple measurement vector (MMV) setup that measures all the snapshots using the same matrix, we consider each snapshot to be measured using a different measurement matrix. We show that this approach reduces the total number of measurements, especially when the rank of the matrix is much smaller than than its sparsity. Our experiments in the context of dynamic imaging shows that this approach is very useful in realizing free breathing cardiac MRI.} \\
	

\section{Introduction} 
The recovery of matrices that are simultaneously low-rank and jointly sparse from its few measurements have received a lot of attention in the recent years. \cite{ktPCA,lingala2012accelerated,liang}. Our interest in this problem is motivated by its application in cardiovascular dynamic MRI, where the image frames are highly correlated. The cardiac signal can be modeled as a low rank matrix $\mathbf X$ with each column as the list of $n$ voxels and each row as the voxel time profile over $N$ snapshots.  Similarly the non-zero locations of the matrix in appropriate transform/frame domains (e.g. wavelet, gradient) are roughly the same in different frame \cite{zhao2011further}. Joint sparsity is a special case of group sparsity where only some groups have non zero coefficients.  For a slowly varying temporal profile of the transformed data we also have joint sparsity.
Our experiments show that the use of joint sparsity and low-rank priors (k-t SLR) provide a considerable reduction in the number of measurements required to recover the matrix.  While there is considerable theoretical progress in problems such as recovering jointly sparse vectors or low-rank matrices, the recovery of matrices that are simultaneously low-rank and jointly sparse have received considerably less attention.

	In \cite{golbabaee2012compressed}, Golbabee et al., have recently developed theoretical guarantees for the recovery of a matrix of rank $r$ and has only $k$ non-zero columns using low rank and joint sparsity priors from its random Gaussian dense measurements. Unfortunately, the dense measurement scheme, where each measurement is a linear combination of all matrix entries is not practical in dynamic imaging; each measurement can only depend on a single column of the matrix.  In contrast, the classical  MMV (multiple measurement vector) problem deals with the recovery of $\mathbf X^{n \times N}$ (N snapshots of length n) from its undersampled measurements $\mathbf Y=\mathbf A\mathbf X$. The matrix $\mathbf X$ is assumed to have a rank of $r$, while the columns are assumed to be jointly sparse with a sparsity of $k$. The spark requirement on the measurement matrix $\mathbf A$ to uniquely identify $\mathbf X$ is ${\rm spark}(\mathbf A)>2k-r+1$ \cite{chen2006theoretical}. While this condition is more relaxed than the identifiability conditions in the single measurement vector (SMV) setting (${\rm spark}(\mathbf A)>2k$), the gains over SMV are marginal when $r<<k$. This is not desirable since the columns of the matrix to be highly correlated and hence more compact; one would expect higher gains in this case. The main problem with this setting is the lack of diversity in the measurements. 
	
	We introduce a hybrid acquisition scheme, where each column $\mathbf x_i$ is measured using the measurement matrix $\mathbf A_i$. All $\mathbf A_i$'s have a common submatrix $\boldsymbol\Phi$, while the remaining rows of $\mathbf A_i$ specified by $\mathbf B_i$ are different. We also introduce a novel two-step reconstruction strategy. In the first step, we estimate the row-subspace using the common measurements $\boldsymbol\Phi$. Once the subspace is estimated, the subspace aware recovery of the matrix simplifies to a least square problem. This approach is considerably simpler than convex optimization schemes using simultaneously structured priors. We introduce perfect recovery guarantees under simple conditions on the measurement matrices and the signal matrix. Our results show that perfect recovery of the matrix is possible when the number of measurements are of the order of the degrees of freedom in the low-rank and sparse matrix. Considering that the dense measurement scheme is not feasible in dynamic imaging applications, this work is quite significant.
	
		The proposed framework to recovery of matrices that are simultaneously jointly-sparse and low-rank is ideally suited to accelerate free breathing cardiac CINE MRI data. The joint sparsity assumption implies that all the snapshots share the same sparsity pattern and thus the the data matrix has only some significant rows and other rows are negligible. Very often, the rank $r$ is much smaller than the joint sparsity $k$. The slow nature of MR acquisitions restricts the number of measurements that can be acquired during each time frame. At the same time, the acquisition is flexible enough to allow the use of different measurement matrices to measure each time frame, without any additional penalty. We propose to capitalize this flexibility to considerably accelerate the acquisition. \textcolor{black}{Specifically, we use the sampling pattern shown in Fig. 1 to acquire the data}.
	\begin{figure}
	\centering
     {\includegraphics[width=.4\textwidth,height = 3.5 cm]{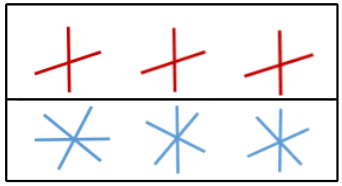}}
\caption{Common (top) and variable (bottom) sampling lines} 
\label{x1rec}
\end{figure}
	\section{Proposed Approach} 
	We propose an acquisition scheme, where the $i^{\rm th}$ frame is measured as 
	\begin{equation}
	\label{maineqns}
	\underbrace{\left[ \begin{array}{c} \mathbf z_{i} \\ \mathbf y_i   \end{array} \right ]}_{ \overline{y}_i}=\underbrace{\left [ \begin{array}{c} \boldsymbol \Phi \\ \mathbf A_i    \end{array} \right ]}_{\mathbf A_{i}}~\mathbf x_{i}.
	\end{equation}
	Note that the measurement matrix $\mathbf \Phi$ is common for all snapshots. We propose to estimate the right subspace of $\mathbf X$ from $\mathbf Z^{H}~\mathbf Z $ , where $Z$ is obtained from $z_i$. Once the right subspace of $\mathbf X$ is estimated, we propose to perform a right subspace aware recovery of $\mathbf X$ from the measurements $\mathbf y_{i}=\mathbf A_{i}~\mathbf x_{i}; ~i=0,..,N-1$. These are the variable measurements and they differ across frames.
	\subsection{Recovery of the right singular matrix}
	In \cite{sketched} the authors look into the preservation of the right singular vectors of a matrix after it is treated by an undersampled matrix. The arrive at a condition for a sensing matrix that satisfies the JL lemma. If the rows of that sensing matrix are of the order of the rank of the matrix for the singular values and vectors to be approximately preserved with a certain failure probability. We look for a similar condition when the sparsity is also coming into picture. 
    
	 Consider jointly sparse matrix $\mathbf X$ and its truncated singular value decomposition, $ \mbf U \mbf \Sigma \mbf V $ where $\bs \Sigma$ is diagonal, non singular matrix in $\mathbb{R}^{r \times r}$; $~ \mbf U \in \mathbb R ^{n \times r}, n>r, \mbf U \mbf U^H=\mbf I;$  $\mbf V \in \mathbb R ^{r \times N}, N>r and, \mbf V \mbf V^H=\mbf I.$ \\ \\ Consider $\mbf Z=\bs \Phi \mbf X = \bs \Phi \mbf U \bs \Sigma \mbf V $ \\ \\Now, $\mbf Z^H \mbf Z=\mbf V^H \bs \Sigma^H \mbf U^H \bs \Phi^H \bs \Phi \mbf U \bs \Sigma \mbf V = \mbf V^H \mbf B^H \mbf B \mbf V,$ where $\mbf B=\bs \Phi \mbf U \bs \Sigma.$
	\noindent
	\begin{thm}
	Every square root of $\mbf Z^H \mbf Z$ can be expressed as $\mbf R \mbf V$ where $\mbf R \in \mathbb R ^{r \times r}$ . Further if $spark(\bs \Phi) \geq  k + 1$ then $\mbf R$ is non zero. If $spark(\bs \Phi)< k$ then $\bs \exists$ jointly sparse $\mbf X$ with rank $r$ s.t $det(\mbf R) = 0.$
	\end{thm} 
	 If $\mbf B$ is full rank, $\mbf B^H \mbf B$ is positive definite and has a unique symmetric square root.
	
	
	
	 Thus for $\mbf Z^H \mbf Z=\mbf V^H \mbf B^H \mbf B \mbf V=\mbf Q^H \mbf Q$ \\ where $\mbf Q=\mbf R \mbf V$ and $\mbf B^H \mbf B$ has a unique square root $\mbf R$ then the right singular matrix of $\mbf V$ can be obtained to a factor $\mbf R$ using $\mbf Z$ iff $\bs \Phi$ has a spark $\geq  k+1.$ \\

	From our experiments we found that $spark(\bs \Phi \geq k+1)$ is strong in practice.  The projection error between original and estimated subspace is falling to 0 when the number of measurements or rows in $\bs \Phi = r.$ See Fig 2 for illustration.
	\noindent Noting this observation we devise the following, the proof of which comes from a determinant based argument provided in \cite{detpolynomialroots}.
	\begin{thm}
	The right singular matrix $\mathbf V$ of any specific matrix $\mathbf X$ can be uniquely recovered from the measurements $\mathbf Z = \boldsymbol \Phi \mathbf X$ for almost all matrices $\mathbf \Phi \in \mathbb C^{r \times n}$.
	\end{thm} 
	
	

			\noindent \textbf{Theorem 1} is a worst case condition that guarantees the recovery for any $ k$ jointly sparse $\mathbf X$, while \textbf{Theorem 2} is sufficient for almost all matrices $\mathbf X$. 
	
\label{x1rec}

\begin{figure}
     {\includegraphics[width=.5\textwidth]{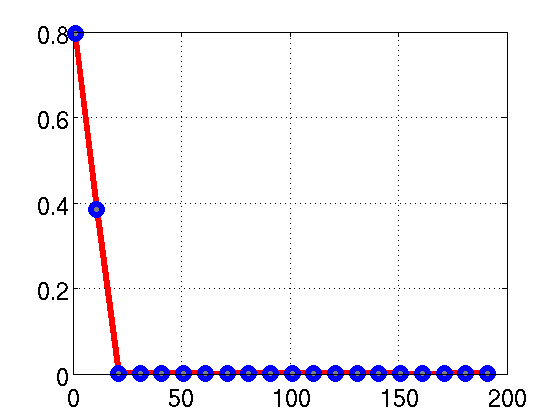}}
\caption{Subspace estimation error vs \# common measurements on Pincat data of rank 20.} 
\label{x1rec}
\end{figure}
	
	\subsection{Subspace aware recovery of $\mathbf X$}
	After the basis vectors of the right subspace specified by the columns of $\mathbf Q \in \mathbb C^{r \times N}$ are estimated, we can express the matrix $\mathbf X$ as $\mathbf X = \mathbf P\mathbf Q$, where $\mathbf P \in  \mathbb C^{n \times r}$ is an appropriate coefficient matrix. 
    From (1) we have, $y_i = \mbf A_i \mbf P q_i$ and for the 1st observation,
    $$y_1 = \mbf A_1 \mbf P q_{11} + \mbf A_1 \mbf P q_{12} + \cdots + \mbf A_1 \mbf P q_{1r} $$
    \begin{equation}
    y_1 = [q_{11} \mbf A_1 ~ q_{12} \mbf A_1 ~ \cdots ~ q_{1r} \mbf A_1  ]
    \left[
\begin{array}{c}
\mathbf p_1\\
\mathbf p_2\\
\vdots\\
\mathbf p_r\\
\end{array}
\right] \nonumber
    \end{equation}
    Vectorizing all the $y_i$ s we obtain,
\begin{eqnarray}
\label{syseq}
	\underbrace{\left [ \begin{array}{c} \mathbf y_1  \\\vdots   \\ \mathbf y_N \end{array} \right]}_{{\rm vec}(\mathbf Y)}
	= \underbrace{\left [\begin{array}{ccc} q_{11}~\mathbf A_1 &\cdots& q_{r1}~\mathbf A_1 \\
	\vdots \\ q_{1N}~\mathbf A_N &\cdots& q_{rN}~\mathbf A_N \end{array} \right]}_{\mathbf B}
	\underbrace{\left [\begin{array}{c}\ \mathbf p_1 \\\vdots   \\ \mathbf p_r \end{array} \right]}_{{\rm vec}(\mathbf P)}
	\end{eqnarray}
    
	Note that the sparsity of ${\rm vec}(\mathbf P)$ is $k r$. Thus, ${\rm vec}(\mathbf P)$ can be uniquely identified if ${\rm spark}(\mathbf B) = 2kr$. We conjecture that this is in fact pessimistic and unique recovery can be achieved if ${\rm spark}(\mathbf A_{i})\geq 2kr/N$. \\
\begin{thm}	
	Assume that $N=(p+1)r$, where $p$ is an arbitrary integer. Let
\begin{eqnarray}\nonumber
\mathbf C_1 & = &\mathbf A_1 = \mathbf A_2  .. = \mathbf A_r\\\nonumber
 \mathbf C_2 & = & \mathbf A_{r+1} = \mathbf A_{r+2}  .. = \mathbf A_{2r}\\\nonumber
\vdots\\\label{clustering}
\mathbf C_{p} &=& \mathbf A_{pr+1} = \mathbf A_{pr+2}  .. = \mathbf A_{N} 
\end{eqnarray}
Here, $\mathbf C_i \in \mathbb R^{m_i \times n}; i = 1,.. p$. $\mathbf P$ can be uniquely determined from (2) if spark($\mathbf X$) = $r + 1$ and
\begin{equation}
{\rm spark}\left(\left[
\begin{array}{c}
\mathbf C_1\\
\mathbf C_2\\
\vdots\\
\mathbf C_p\\
\end{array}
\right]
\right) \geq 2k.
\end{equation}
\end{thm}
We present here a proof for the rank one case. The proof for a general low rank $r$ case is omitted here for brevity. 

\begin{proof}
Consider $r =1.$ 
$y_i = \mbf A_i x_i = \mbf A_i u_1 v_i$, where $u$ has only one column as the rank is 1 and $v_i$s are scalars. We have, $y_i/v_i = \mbf A_i u_1.$ From the compressed sensing results on  unique recovery of a single unknown sparse vector we need, 
\begin{equation}
{\rm spark}\left(\left[
\begin{array}{c}
\mathbf A_1\\
\mathbf A_2\\
\vdots\\
\mathbf A_N\\
\end{array}
\right]
\right) \geq 2k.
\end{equation}

\end{proof}

Note that $\mathbf X$ and $\mathbf Q$  share the same null space. So from the spark condition $spark(\mbf X) = r+1 $ in Theorem 3, we get that any $r$ columns of $\mbf Q$ forms an non singular matrix. If we consider the clustering in (3), we only need the $\mbf Q$ submatrix formed by picking the rows corresponding to the indices in a group in (3) to be non singular. Simply, we want consecutive $r$ rows of $\mbf Q^H$ to be full rank. This requirement comes from the proof of the general case.

This is only required for the sufficient condition. In practice as shown by our simulation results we do not need that spark condition on $\mathbf X$. It is just to attain the theoretical claim.
	\subsection{Effect of sampling pattern }
The clustering specified in (3) assumes that every consecutive r frames are observed by same set of variable measurements. We have devised this clustering of sensing matrices to arrive at the sufficient condition in (4). As specified above, our goal is to ensure that the frames which are operated by similar sensing matrices or in other words the frames that are operated by matrices in the same group should form an invertible matrix. The Q submatrix such formed should be non singular. The reconstruction of simulated data satisfying (3) is shown in Fig 3. We see perfect recovery in the noiseless case on the left. 

This clustering can also be done differently. Possibilities include repeating the $\mathbf A_i s$ after every r frame hence generating a periodicity in the observation matrix. Another possibility is taking the given cluster and permuting it across the N frames i.e., randomly shuffling the set of sensing matrices such that a certain $r$ of them are same. We want to look at how the reconstructions look like if we violate from the pattern we have in (3), i.e., for a given number of measurements how the three patterns perform. This depends on the knowledge of the structure of the data we are trying to reconstruct.

We will highlight this effect in Fig 4.  The chosen toy Pincat dataset is periodic in nature with a rank 20. So every $20^{th}$ frame belongs to the same cardiac phase. So the matrix formed by stacking every $20^{th}$ frame is rank 1. Clearly, the clustering in (3) is favorable here. From the possibilities mentioned above, the periodic clustering will fail as the $\mathbf Q$ submatrix formed using every $r = 20^{th}$ frame will be rank deficient or in this case rank 1.  Random permutation works a bit better than the periodic clustering as it brings some complimentary information. The rank of the $\mbf Q$ submatrix formed by the random clustering is still higher than 1 but not full rank, hence the recovery isn't perfect.

So the underlying fact remains that we can choose any clustering we want provided the $C_i s$ formed satisfy (4) and we have the idea of roughly how the time frames behave (periodicity) for a particular dataset. And lastly, this is only to comply with the theory. In practice as we will see in Fig. 5 , even without satisfying this spark condition $spark(\mbf X) = r+1$ we can get a good recovery.

\begin{figure}[t!]
     \subfigure[{Noiseless}]{\includegraphics[width=.24\textwidth,height = 4 cm]{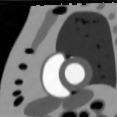}}
     \subfigure[{SNR : 35 dB}]{\includegraphics[width=.24\textwidth,height = 4 cm]{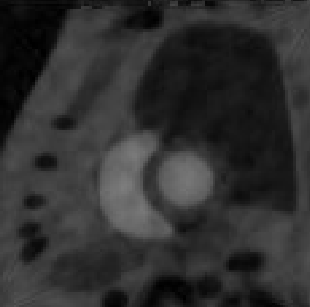}}
\caption{Reconstructing Pincat data using our two step strategy. Left: in a noiseless setting. Right: with a SNR of 35 dB } 
\label{x1rec}
\end{figure}

\begin{figure}[t!]
     \subfigure[{Periodic: Every $r^{th}$}]{\includegraphics[width=.24\textwidth,height = 4 cm]{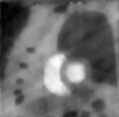}}
     \subfigure[{Permuting the groups across all the frames}]{\includegraphics[width=.24\textwidth,height = 4 cm]{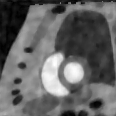}}
\caption{Effect of clustering the similar sensing matrices on the reconstruction:  placing them periodically (left) and permuting them randomly (right)} 
\label{x1rec}
\end{figure}

		\subsection{Measurements required}
		For a matrix of dimension $n \times N$ and rank r, classical MMV scheme requires  $s \geq (2k-r+1)N$ measurements for its unique recovery. Dense measurement scheme requires $s$ of the order of the degrees of freedom in a matrix. Variable and common measurements in the proposed scheme requires $s \geq (2k-r+N)r$ for unique recovery. Combining the above, $\mathbf X$ can be uniquely identified if the average number of measurements per snapshot is of the order of $s=r + 2kr/N$. Note that $s \rightarrow r$ as $N\rightarrow \infty$.  This is a quite considerable gain over the number required for MMV identifiability, specified by $2k-r+1$, especially when $r<<k$.

			
		\section{Algorithm}
		We use Alternating direction method of multipliers (ADMM) \cite{afonso2011augmented,hestenes} which combines the augmented Lagrangian method and variable splitting techniques to solve the constrained optimization problem. The problem is posed as:
\begin{eqnarray}
\underset{\mathbf P}{\operatorname{\mathbf argmin}} ~ ||\mathbf B~ \mathbf vec(\mathbf P) -vec(\mathbf Y)||_2^2 + \nonumber  \\
\boldsymbol \lambda \sum_{xy} \sqrt{\sum_t {\mathbf P_x(t)^2+\mathbf P_y(t)^2} }
\end{eqnarray}	 With the V's supplied by the 1st step we solve for P's by creating an auxiliary variable G. After the variable splitting we have the Lagrangian,	  

\begin{eqnarray}
 \mathcal{L} =  ||\mathbf B ~vec(\mathbf P) -vec(\mathbf Y)||_2^2 + \nonumber  \\ \boldsymbol \lambda \sum_{xy} \sqrt{\sum_t {\mathbf G_x(t)^2 +\mathbf G_y(t)^2} }\nonumber  \\
   +\boldsymbol \lambda \frac{\boldsymbol \beta}{2} ||\nabla_x \mathbf P-\mathbf G_x||_2^2 + \boldsymbol \lambda \frac{\boldsymbol \beta}{2} ||\nabla_y \mathbf P-\mathbf G_y||_2^2 \nonumber  \\
   + <\Lambda_x, \mathbf G_x-\nabla_x \mathbf P > + <\Lambda_y, \mathbf G_y-\nabla_y \mathbf P >
\end{eqnarray}
where $\bs \beta$ is the penalty parameter which drives the difference between the auxiliary and original variable to zero, $\bs \lambda$ is the regularization parameter, $\mathbf \Lambda_x , \mathbf \Lambda_y $ are the Lagrange multipliers.
 The P sub-problem is solved using Preconditioned conjugate gradient and G sub-problem is solved by $l_1$ shrinkage. The data is transformed to the gradient domain to make it sparse so the problem is solved with a Total variation (TV) regularization. $l_1$ norm encourages sparsity and $l_2$ norm encourages diversity \cite{sensor}. We enforce joint sparsity by penalizing the $l_1-l_2$ mixed norm of the gradient data.

\begin{figure}
\includegraphics[width=.5\textwidth]{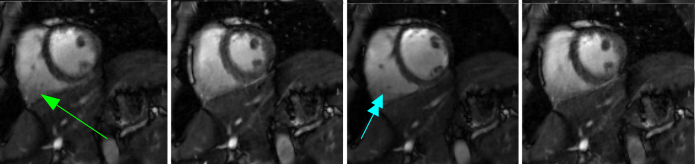}
\caption{Reconstruction of CINE MRI data using our algorithm. Few cardiac phases are shown. The dots pointed out by the arrow shows the same phases reappearing.} 
\label{x1rec}
\end{figure}

\section{Numerical results}
We validate our Theorems on toy Pincat \cite{sharif2007physiologically} phantom data of dimension 128 x 128 x 200 and rank 20. In Fig. 2 we consider a noiseless setting. We plot the projection error between the estimated and actual right subspace vs the number of random samples of the image data in the Fourier domain. Projection error between two subspaces $V_1$ and $V_2$ is defined as:
\begin{equation}
\frac{||(\mathbf I-\mathbf V_1\mathbf V_1^H)\mathbf V_2||_2^2 +||(\mathbf I-\mathbf V_2\mathbf V_2^H)\mathbf V_1||_2^2 }{||\mathbf V_1||_2^2+||\mathbf V_2||_2^2}.
\end{equation} 
We see a drop in the normalized error after the the number of samples equals the rank 20.

In Fig. 3 we look at the subspace aware recovery of the Pincat data using the clustering in (3) which is consistent with the formulation of Theorem 4. This was undersampled with 4 common and 5 variable lines for a noiseless (left) and a noisy setting (SNR of 40 db on the right). We get perfect recovery in the noiseless case. All the Pincat performance comparison can be done with the perfect reconstruction Fig. 3a as reference. We see blurring artifacts in the noisy case Fig. 3b.

Fig. 4 shows the effect of clustering the sampling pattern on the reconstruction. On the left figure the similar sensing matrices are distributed periodically to every $r^{th}$ frame and on the right the clusters are shuffled randomly across all the N frames. We see that for this dataset, following the setting of Theorem 4 better reconstructions are obtained for the same number of measurements as against when those conditions were violated. This indicates that the clustering pattern can be modified based on knowledge of the structure in the data. But after that if (4) is satisfied good recovery is ensured. In all the Pincat results, the first frame is shown from the reconstructed time series.

We apply our algorithm on free breathing CINE data in Fig. 5. A cropped portion of the reconstruction is shown highlighting the heart. The dots indicated by the arrows show the similar cardiac phases. The data was acquired using an SSFP sequence with an $18$ channel coil array, with TR/TE of $4.2/2.1$ ms, matrix size of $512\times512$, FOV of $300$mm$\times300$mm and slice thickness of $5$mm on 3T Siemens Trio scanner.  We considered $12$ radial lines of Fourier space to reconstruct each image frame, $4$ of which were common measurement lines. This corresponds to a temporal resolution of $50$ ms. The acquisition time was $25$ s translates to $500$ image frames. Each radial Fourier line is a 512 pixels long line on the 2D Fourier space data passing through the center of the Fourier space as shown in Fig 1. The blurring is optimal and the cardiac phases are clearly visible. Despite the undersampling our algorithm is providing good reconstruction using the two step strategy. This is significant because MR acquisition is slow and we can't fully sample MR data to get the perfect reconstruction. Thus the strategic undersampling to get good recovery is essential when the data has a low rank and joint sparse structure.

	\section{Conclusion}
	We propose to reconstruct a low rank and sparse matrix from its under sampled measurements. We introduced a two step subspace strategy that recovers the right and left subspaces one after another. Our sufficient conditions show that by using our sampling scheme we are gaining in terms of the number of measurements to uniquely recover the data as compared to the MMV scheme.  Our experiments suggest that considerable undersampling when strategically done, gives good reconstruction in the noiseless case for the Pincat data. In the CINE MRI data where the ground truth is unknown considerable undersampling along with our reconstruction algorithm shows good recovery. The reduced number of required measurements due to undersampling will translate to high accelerations in dynamic MRI. \\

\bibliographystyle{IEEEtran}
\bibliography{refs}

\begin{thebibliography}{10}
\providecommand{\url}[1]{#1}
\csname url@samestyle\endcsname
\providecommand{\newblock}{\relax}
\providecommand{\bibinfo}[2]{#2}
\providecommand{\BIBentrySTDinterwordspacing}{\spaceskip=0pt\relax}
\providecommand{\BIBentryALTinterwordstretchfactor}{4}
\providecommand{\BIBentryALTinterwordspacing}{\spaceskip=\fontdimen2\font plus
\BIBentryALTinterwordstretchfactor\fontdimen3\font minus
  \fontdimen4\font\relax}
\providecommand{\BIBforeignlanguage}[2]{{%
\expandafter\ifx\csname l@#1\endcsname\relax
\typeout{** WARNING: IEEEtran.bst: No hyphenation pattern has been}%
\typeout{** loaded for the language `#1'. Using the pattern for}%
\typeout{** the default language instead.}%
\else
\language=\csname l@#1\endcsname
\fi
#2}}
\providecommand{\BIBdecl}{\relax}
\BIBdecl

\bibitem{ktPCA}
H.~Pedersen, S.~Kozerke, S.~Ringgaard, K.~Nehrke, and W.~Y. Kim, ``{k-t PCA}:
  temporally constrained {k-t BLAST} reconstruction using principal component
  analysis,'' \emph{Magn Reson Med}, vol.~62, no.~3, pp. 706--716, Sep 2009.

\bibitem{lingala2012accelerated}
S.~G. Lingala, E.~DiBella, M.~Jacob \emph{et~al.}, ``Accelerated imaging of
  rest and stress myocardial perfusion mri using multi-coil kt slr: a
  feasibility study,'' \emph{Journal of Cardiovascular Magnetic Resonance},
  vol.~14, no. Suppl 1, p. P239, 2012.

\bibitem{liang}
Z.~Liang, ``Spatiotemporal imaging with partially separable functions,'' in
  \emph{ISBI}, 2007, pp. 181--182.

\bibitem{zhao2011further}
B.~Zhao, J.~Haldar, A.~Christodoulou, and Z.~Liang, ``Further development of
  image reconstruction from highly undersampled (k, t)-space data with joint
  partial separability and sparsity constraints,'' in \emph{Biomedical Imaging:
  From Nano to Macro, 2011 IEEE International Symposium on}.\hskip 1em plus
  0.5em minus 0.4em\relax IEEE, 2011, pp. 1593--1596.

\bibitem{golbabaee2012compressed}
M.~Golbabaee and P.~Vandergheynst, ``Compressed sensing of simultaneous
  low-rank and joint-sparse matrices,'' \emph{arXiv preprint arXiv:1211.5058},
  2012.

\bibitem{chen2006theoretical}
J.~Chen and X.~Huo, ``Theoretical results on sparse representations of
  multiple-measurement vectors,'' \emph{IEEE Transactions on Signal
  Processing}, vol.~54, no.~12, pp. 4634--4643, 2006.

\bibitem{sketched}
A.~C. Gilbert, J.~Y. Park, and M.~B. Wakin, ``Sketched {SVD:} recovering
  spectral features from compressive measurements,'' \emph{Arxiv preprint
  arXiv:1211.0361}, 2012.

\bibitem{detpolynomialroots}
L.~Losaz, ``Singular spaces of matrices and their application in
  combinatorics,'' \emph{Bulletin of the Brazilian Mathmatical Society},
  vol.~20, no.~1, pp. 87--89, 1989.

\bibitem{afonso2011augmented}
M.~Afonso, J.~Bioucas-Dias, and M.~Figueiredo, ``An augmented lagrangian
  approach to the constrained optimization formulation of imaging inverse
  problems,'' \emph{Image Processing, IEEE Transactions on}, vol.~20, no.~3,
  pp. 681--695, 2011.

\bibitem{hestenes}
M.~R. Hestenes, ``Multiplier and gradient methods,'' \emph{Journal of
  Optimization Theory and Applications}, vol.~4, no.~5, pp. 303--320, 1969.

\bibitem{sensor}
J.~Le~Roux, P.~T. Boufounos, K.~Kang, and J.~R. Hershey, ``Source localization
  in reverberant environments using sparse optimization,'' in \emph{Acoustics,
  Speech and Signal Processing (ICASSP), 2013 IEEE International Conference
  on}.\hskip 1em plus 0.5em minus 0.4em\relax IEEE, 2013, pp. 4310 -- 4314.

\bibitem{sharif2007physiologically}
B.~Sharif and Y.~Bresler, ``Physiologically improved ncat phantom (pincat)
  enables in-silico study of the effects of beat-to-beat variability on cardiac
  mr,'' in \emph{Proceedings of the Annual Meeting of ISMRM, Berlin}, 2007, p.
  3418.

\end{thebibliography}
	 
	\end{document}